\documentstyle[11pt, amssymb]{article}

\begin{document}

\begin{flushleft}
{\Large Periods of generalized Tate curves}
\end{flushleft}

\begin{flushleft}
{\large Takashi Ichikawa} 
\end{flushleft}

\begin{flushleft} 
{\it Department of Mathematics, Faculty of Science and Engineering, 
Saga University, Saga 840-8502, Japan} 
\footnote{{\it E-mail address:}  ichikawn@cc.saga-u.ac.jp.} 
\end{flushleft} 

\begin{flushleft}
\end{flushleft}

\begin{flushleft}
\end{flushleft}

\noindent
ABSTRACT 

\noindent
A generalized Tate curve is a universal family of curves 
with fixed genus and degeneration data 
which becomes Schottky uniformized Riemann surfaces and Mumford curves 
by specializing moduli and deformation parameters. 
By considering each generalized Tate curve as a family of degenerating Riemann surfaces, 
we give explicit formulas of the period isomorphism 
between its de Rham and Betti cohomology groups, 
and of the associated objects: 
Gauss-Manin connection, variation of Hodge structure and monodromy weight filtration. 
A remarkable fact is that similar formulas hold also for families of Mumford curves. 
Furthermore, we show that for a generalized Tate with maximally degenerate closed fiber, 
its local unipotent periods can be expressed as power series in the deformation parameters 
whose coefficients are multiple polylogarithm functions. 
This $p$-adic version is also given. 
\vspace{2ex}

\noindent
\hrulefill
\vspace{2ex}

\noindent
{\bf 1. Introduction} 
\vspace{2ex}

The study of period isomorphisms between the de Rham and Betti cohomology groups 
of families of algebraic varieties is one of main subjects in algebraic geometry. 
Their behavior under degeneration of considering varieties is especially interesting, 
and the associated $p$-adic version, unipotent version and $p$-adic unipotent version 
are also studied. 
For a universal family of elliptic curves, 
these objects are described by hypergeometric functions and equations, 
however, there seems no similar description on curves of genus $>1$. 
The aim of this paper is to give such explicit formulas 
for a universal family of general curves which are commonly applicable to 
complex and $p$-adic analytic cases using the theory of generalized Tate curves. 

A generalized Tate curve is introduced in \cite{I1} 
as a higher genus version of the Tate elliptic curve
which gives a universal deformation of degenerate curves of fixed type. 
The base ring of a generalized Tate curve is taken to be primitive as far as possible, 
namely it consists of formal power series in deformation parameters 
whose coefficients belong to the ring over ${\mathbb Z}$ of moduli parameters 
of the degenerate curves. 
By specializing these parameters in ${\mathbb C}$ 
(resp. a nonarchimedean complete valuation field), 
a generalized Tate curve gives families of Schottky uniformized Riemann surfaces 
(resp. Mumford curves). 

In this paper, 
we give explicit formulas of abelian differentials defined on a generalized Tate curve 
called {\it universal differentials}, 
and show their variational formulas which imply an important fact 
that the universal differentials are stable, 
namely they have only logarithmic poles on the special fiber. 
By regarding the generalized Tate curve as a family of degenerating Riemann surfaces, 
we give the explicit description of the period isomorphism 
between the de Rham and Betti cohomology groups, 
and of the associated Gauss-Manin connection, 
variation of Hodge structure and monodromy weight filtration. 
This description is written by multiple moduli and deformation parameters 
in the category of formal geometry over ${\mathbb Q}$, 
and is same to formulas for Mumford curves 
which is studied by Gerritzen \cite{G} and de Shalit \cite{dS1, dS2} 
using Coleman integration \cite{C, CdS} in the special case of $p$-adic base fields. 

As a further application of the above result,  
we study unipotent periods which give isomorphisms between 
the de Rham and Betti unipotent fundamental groups of curves. 
Using the stability of universal differentials, 
for a generalized Tate curve with maximally degenerate closed fiber, 
its local unipotent periods can be expressed as power series in the deformation parameters
whose coefficients are multiple polylogarithm functions. 
This result extends results of Deligne \cite{D2}, Drinfel'd \cite{Dr} and Hain \cite{H2} 
on the Knizhnik-Zamolodchikov-Bernard equation in the case of genus $\leq 1$. 
Furthermore, 
a similar result on multiple zeta values is obtained by Brown \cite{Br} and Banks-Panzer-Pym \cite{BPP} 
in the genus $0$ case, 
and is applied to a conjecture of Kontsevich \cite{K} on deformation quantization. 
Similar descriptions are given for local $p$-adic unipotent periods 
based on result of Furusho \cite{F1, F2} which give examples of Coleman integrals in families 
considered by Besser \cite{Be}. 
Therefore, we may give a hope which states, roughly speaking, 
that if a family of curves with any genus degenerates maximally, 
then their motivic structure has a perturbative expansion in the deformation parameters 
whose coefficients are described by mixed Tate motives over ${\mathbb Z}$. 
\vspace{4ex}

\noindent
{\bf 2. Generalized Tate curve}  
\vspace{2ex}

\noindent
{\it 2.1. Schottky uniformization} 
\vspace{2ex}

A Schottky group $\Gamma$ of rank $g$ is a free group with generators 
$\gamma_{i} \in PGL_{2}({\mathbb C})$ $(1 \leq i \leq g)$ which map Jordan curves  
$C_{i} \subset {\mathbb P}^{1}_{\mathbb C} = {\mathbb C} \cup \{ \infty \}$ 
to other Jordan curves $C_{-i} \subset {\mathbb P}^{1}_{\mathbb C}$ 
with orientation reversed, 
where $C_{\pm 1},..., C_{\pm g}$ with their interiors are mutually disjoint. 
Each element $\gamma \in \Gamma - \{ 1 \}$ is conjugate 
to an element of $PGL_{2}({\mathbb C})$ sending $z$ to $\beta_{\gamma} z$ 
for some $\beta_{\gamma} \in {\mathbb C}^{\times}$ with $|\beta_{\gamma}| < 1$ 
which is called the {\it multiplier} of $\gamma$. 
Therefore, one has 
$$
\frac{\gamma(z) - \alpha_{\gamma}}{z - \alpha_{\gamma}} = 
\beta_{\gamma} \frac{\gamma(z) - \alpha'_{\gamma}}{z - \alpha'_{\gamma}} 
$$
for some element $\alpha_{\gamma}, \alpha'_{\gamma}$ of 
${\mathbb P}^{1}_{\mathbb C}$ 
called the {\it attractive}, {\it repulsive} fixed points of $\gamma$ 
respectively. 
Then the discontinuity set 
$\Omega_{\Gamma} \subset {\mathbb P}^{1}_{\mathbb C}$ 
under the action of $\Gamma$ has a fundamental domain $D_{\Gamma}$ 
which is given by the complement of the union of the interiors of $C_{\pm 1},..., C_{\pm g}$. 
The quotient space $R_{\Gamma} = \Omega_{\Gamma}/\Gamma$ 
is a (compact) Riemann surface of genus $g$ 
which is called Schottky uniformized by $\Gamma$ (cf. \cite{Scho}). 
Furthermore, by a result of Koebe, 
every Riemann surface of genus $g$ can be represented in this manner. 
\vspace{2ex}

\noindent
{\it 2.2. Degenerate curve} 
\vspace{2ex}

We review a correspondence 
between certain graphs and degenerate pointed curves, 
where a (pointed) curve is called {\it degenerate} if it is a stable (pointed) curve and 
the normalization of its irreducible components are all projective (pointed) lines. 
A {\it graph} $\Delta = (V, E, T)$ means a collection  
of 3 finite sets $V$ of vertices, $E$ of edges, $T$ of tails 
and 2 boundary maps 
$$
b : T \rightarrow V, 
\ \ b : E \longrightarrow \left( V \cup \{ \mbox{unordered pairs of elements of $V$} \} \right) 
$$
such that the geometric realization of $\Delta$ is connected. 
A graph $\Delta$ is called {\it stable} 
if its each vertex has at least $3$ branches. 
Then for a degenerate pointed curve, 
its dual graph $\Delta = (V, E, T)$ is given by the correspondence: 
$$
\begin{array}{lcl}
V & \longleftrightarrow & 
\{ \mbox{irreducible components of the curve} \}, 
\\
E & \longleftrightarrow & 
\{ \mbox{singular points on the curve} \}, 
\\
T & \longleftrightarrow & 
\{ \mbox{marked points on the curve} \} 
\end{array}
$$
such that an edge (resp. a tail) of $\Delta$ has a vertex as its boundary 
if the corresponding singular (resp. marked) point belongs 
to the corresponding component. 
Denote by $\sharp X$ the number of elements of a finite set $X$. 
Under fixing a bijection 
$\nu : T \stackrel{\sim}{\rightarrow} \{ 1, ... , \sharp T \}$, 
which we call a numbering of $T$, 
a stable graph $\Delta = (V, E, T)$ becomes the dual graph 
of a degenerate $\sharp T$-pointed curve of genus 
${\rm rank}_{\mathbb Z} H_{1}(\Delta, {\mathbb Z})$ 
and that each tail $h \in T$ corresponds to the $\nu(h)$th marked point. 
In particular, a stable graph without tail is the dual graph of 
a degenerate (unpointed) curve by this correspondence. 
If $\Delta$ is trivalent, i.e. any vertex of $\Delta$ has just $3$ branches, 
then a degenerate $\sharp T$-pointed curve with dual graph $\Delta$ 
is maximally degenerate. 

An {\it orientation} of $\Delta = (V, E, T)$ means 
giving an orientation of each $e \in E$. 
Under an orientation of $\Delta$, 
denote by 
$$
\pm E 
\ = \ 
\{ e, -e \ | \ e \in E \} 
$$
the set of oriented edges, 
and by $v_{h}$ the terminal vertex of $h \in \pm E$ (resp. the boundary vertex of $h \in T)$. 
For each $h \in \pm E$, 
denote by let $| h | \in E$ be the edge $h$ without orientation. 
\vspace{2ex}

\noindent
{\it 2.3. Generalized Tate curve} 
\vspace{2ex}

Let $\Delta = (V, E, T)$ be a stable graph. 
Fix an orientation of $\Delta$, 
and take a subset ${\cal E}$ of $\pm E \cup T$ 
whose complement ${\cal E}_{\infty}$ satisfies the condition that 
$$
\pm E \cap {\cal E}_{\infty} \cap \{ -h \ | \ h \in {\cal E}_{\infty} \} 
\ = \ 
\emptyset, 
$$ 
and that $v_{h} \neq v_{h'}$ for any distinct $h, h' \in {\cal E}_{\infty}$. 
We attach variables $x_{h}$ for $h \in {\cal E}$ and $y_{e} = y_{-e}$ for $e \in E$. 
Let $A_{0}$ be the ${\mathbb Z}$-algebra generated by $x_{h}$ $(h \in {\cal E})$, 
$1/(x_{e} - x_{-e})$ $(e, -e \in {\cal E})$ and $1/(x_{h} - x_{h'})$ 
$(h, h' \in {\cal E}$ with $h \neq h'$ and $v_{h} = v_{h'})$, 
and let 
$$ 
A_{\Delta} \ = \ A_{0} [[y_{e} \ (e \in E)]], \ \ 
B_{\Delta} \ = \ A_{\Delta} \left[ \prod_{e \in E} y_{e}^{-1} \right]. 
$$ 
According to \cite[Section 2]{I1}, 
we construct the universal Schottky group $\Gamma$ 
associated with oriented $\Delta$ and ${\cal E}$ as follows. 
For $h \in \pm E$, put 
\begin{eqnarray*}
\phi_{h} 
& = & 
\left( \begin{array}{cc} x_{h} & x_{-h} \\ 1 & 1 \end{array} \right) 
\left( \begin{array}{cc} 1 & 0 \\ 0 & y_{h} \end{array} \right) 
\left( \begin{array}{cc} x_{h} & x_{-h} \\ 1 & 1 \end{array} \right)^{-1} 
\\ 
& = & 
\frac{1}{x_{h} - x_{-h}} 
\left\{ \left( \begin{array}{cc} x_{h} & - x_{h} x_{-h} \\ 1 & - x_{-h} \end{array} \right) - 
\left( \begin{array}{cc} x_{-h} & - x_{h} x_{-h} \\ 1 & - x_{h} \end{array} \right) y_{h} \right\}, 
\end{eqnarray*} 
where $x_{h}$ (resp. $x_{-h})$ means $\infty$ 
if $h$ (resp. $-h)$ belongs to ${\cal E}_{\infty}$. 
This gives an element of $PGL_{2}(B_{\Delta}) = GL_{2}(B_{\Delta})/B_{\Delta}^{\times}$ 
which we denote by the same symbol, 
and satisfies
$$
\frac{\phi_{h}(z) - x_{h}}{z - x_{h}} 
\ = \ 
y_{h} \frac{\phi_{h}(z) - x_{-h}}{z - x_{-h}} 
\ \ (z \in {\mathbb P}^{1}), 
$$
where $PGL_{2}$ acts on ${\mathbb P}^{1}$ by linear fractional transformation. 
\vspace{2ex}

\noindent
{\bf Proposition 2.1.} 
\begin{it}
Let $\phi$ be a product $\phi_{h(1)} \cdots \phi_{h(l)}$ with 
$v_{-h(i)} = v_{h(i+1)}$ $(1 \leq i \leq l-1)$ 
which is reduced in the sense that $h(i) \neq -h(i+1)$ $(1 \leq i \leq l-1)$, 
and put $y_{\phi} = y_{h(1)} \cdots y_{h(l)}$. 

{\rm (1)} 
One has 
$\displaystyle \phi(z) - x_{h(1)} \in y_{h(1)} 
\left( A_{0} \left[ z, \prod_{h \in \pm E} (z - x_{h})^{-1} \right] [[ y_{e} \ (e \in E) ]] \right)$. 

{\rm (2)} 
If $a \in A_{\Delta}$ satisfies $a - x_{-h(l)} \in A_{\Delta}^{\times}$, 
then $\phi(a) - x_{h(1)} \in I$. 
Furthermore, if $a' - x_{-h(l)} \in A_{\Delta}^{\times}$, 
then $\phi(a) - \phi(a') \in (a - a') y_{\phi} A_{\Delta}^{\times}$. 

{\rm (3)} 
One has 
$\displaystyle 
\frac{d \phi(z)}{dz} \in y_{\phi} 
\left( A_{0} \left[ \prod_{h \in \pm E} (z - x_{h})^{-1} \right] [[ y_{e} \ (e \in E) ]] \right)$. 
\end{it} 
\vspace{2ex}

\noindent
{\bf Proof.}
Since the assertion (2) is proved in \cite[Lemma 1.2]{I1}, 
we will prove (1) and (3). 
Put 
$$
\phi = \left( \begin{array}{cc} a_{\phi} & b_{\phi} \\ c_{\phi} & d_{\phi} \end{array} \right). 
$$
Since 
$$
\left( \begin{array}{cc} \alpha & - \alpha \beta \\ 1 & - \beta \end{array} \right)
\left( \begin{array}{cc} \gamma & - \gamma \delta \\ 1 & - \delta \end{array} \right) 
= (\gamma - \beta) 
\left( \begin{array}{cc} \alpha & - \alpha \beta \\ 1 & - \delta \end{array} \right), 
$$
$a_{\phi}$, $b_{\phi}$, $c_{\phi}$ and $d_{\phi}$ are elements of $A_{\Delta}$ 
whose constant terms are $x_{h(1)} t$, $- x_{h(1)} x_{-h(l)} t$, $t$ and $- x_{-h(l)} t$ 
respectively,   
where 
$$
t = \frac{\prod_{s = 2}^{l} \left( x_{h(s)} - x_{-h(s-1)} \right)}
{\prod_{s = 1}^{l} \left( x_{h(s)} - x_{-h(s)} \right)} \in A_{\Delta}^{\times}. 
$$
Then $c_{\phi} z + d_{\phi} = t (z - x_{-h(l)}) + \cdots$, 
and hence  
$$
\phi(z) - x_{h(1)} \in 
A_{0} \left[ z, \prod_{h \in \pm E} (z - x_{h})^{-1} \right] [[ y_{e} \ (e \in E) ]]. 
$$
In order to prove (1), we may assume that $l = 1$, 
and then $\phi(z) = x_{h(1)}$ under $y_{h(1)} = 0$. 
Therefore, the assertion (1) holds.  
The assertion (3) follows from 
$$
\frac{d \phi(z)}{dz} = \frac{\det(\phi)}{(c_{\phi} z + d_{\phi})^{2}} 
= \frac{\prod_{s = 1}^{l} \det(\phi_{h(s)})}{(c_{\phi} z + d_{\phi})^{2}} 
= \frac{\prod_{s = 1}^{l} y_{h(s)}}{(c_{\phi} z + d_{\phi})^{2}}, 
$$
and the above calculation. 
\ $\square$ 
\vspace{2ex}

For any reduced path $\rho = h(1) \cdot h(2) \cdots h(l)$ 
which is the product of oriented edges $h(1), ... ,h(l)$ such that $v_{h(i)} = v_{-h(i+1)}$, 
one can associate an element $\rho^{*}$ of $PGL_{2}(B_{\Delta})$ 
having reduced expression 
$\phi_{h(l)} \phi_{h(l-1)} \cdots \phi_{h(1)}$. 
Fix a base vertex $v_{b}$ of $V$, 
and consider the fundamental group 
$\pi_{1} (\Delta, v_{b})$ which is a free group 
of rank $g = {\rm rank}_{\mathbb Z} H_{1}(\Delta, {\mathbb Z})$. 
Then the correspondence $\rho \mapsto \rho^{*}$ 
gives an injective anti-homomorphism 
$\pi_{1} (\Delta, v_{b}) \rightarrow PGL_{2}(B_{\Delta})$ 
whose image is denoted by $\Gamma_{\Delta}$. 

It is shown in \cite[Section 3]{I1} and \cite[1.4]{I2} 
(see also \cite[Section 2]{IhN} when $\Delta$ is trivalent and has no loop) 
that for any stable graph $\Delta$, 
there exists a stable pointed curve $C_{\Delta}$ of genus $g$ over $A_{\Delta}$ 
which satisfies the following: 

\begin{itemize}

\item 
The closed fiber $C_{\Delta} \otimes_{A_{\Delta}} A_{0}$ of $C_{\Delta}$ 
obtained by substituting $y_{e} = 0$ $(e \in E)$ 
becomes the degenerate pointed curve over $A_{0}$ with dual graph $\Delta$ which is 
obtained from the collection of $P_{v} := {\mathbb P}^{1}_{A_{0}}$ $(v \in V)$ 
by identifying the points $x_{e} \in P_{v_{e}}$ and $x_{-e} \in P_{v_{-e}}$ ($e \in E$), 
where $x_{h}$ denotes $\infty$ if $h \in {\cal E}_{\infty}$. 

\item 
$C_{\Delta}$ gives rise to a universal deformation 
of degenerate pointed curves with dual graph $\Delta$. 
More precisely, if $R$ is a noetherian and normal complete local ring with residue field $k$, 
and $C$ is a stable pointed curve over $R$ with nonsingular generic fiber 
such that the closed fiber $C \otimes_{R} k$ is a degenerate pointed curve 
with dual graph $\Delta$, in which all double points are $k$-rational, 
then there exists a ring homomorphism $A_{\Delta} \rightarrow R$ 
giving $C_{\Delta} \otimes_{A_{\Delta}} R \cong C$.  

\item 
$C_{\Delta} \otimes_{A_{\Delta}} B_{\Delta}$ is smooth over $B_{\Delta}$ 
and is Mumford uniformized (cf. \cite{Mu}) by $\Gamma$. 

\item 
Take $x_{h}$ $(h \in {\cal E})$ as complex numbers such that $x_{e} \neq x_{-e}$ 
and that $x_{h} \neq x_{h'}$ if $h \neq h'$ and $v_{h} = v_{h'}$, 
and take $y_{e}$ $(e \in E)$ as sufficiently small nonzero complex numbers. 
Then $C_{\Delta}$ becomes a pointed Riemann surface which is 
Schottky uniformized by the Schottky group $\Gamma$ over ${\mathbb C}$ 
obtained from $\Gamma_{\Delta}$. 

\end{itemize}
We review the construction of $C_{\Delta}$ given in \cite[Theorem 3.5]{I1}. 
Let $T_{\Delta}$ be the tree obtained as the universal cover of $\Delta$, 
and denote by ${\cal P}_{T_{\Delta}}$ be the formal scheme 
as the union of ${\mathbb P}^{1}_{A_{\Delta}}$'s indexed by vertices of $T_{\Delta}$ 
under the $B_{\Delta}$-isomorphism by $\phi_{e}$ $(e \in E)$. 
Then it is shown in \cite[Theorem 3.5]{I1} that 
$C_{\Delta}$ is the formal scheme theoretic quotient of ${\cal P}_{T_{\Delta}}$ by $\Gamma_{\Delta}$. 
\vspace{4ex}

\noindent 
{\bf 3. Universal differential} 
\vspace{2ex}

\noindent 
{\it 3.1. Universal differential} 
\vspace{2ex}

We define universal differentials on a generalized Tate curve. 
Let $\Delta$ be a stable graph, 
and 
let 
$\Gamma_{\Delta} = {\rm Im} \left( \pi_{1}(\Delta, v_{b}) \rightarrow 
PGL_{2}(B_{\Delta}) \right)$ 
be the universal Schottky group as above. 
Then it is shown in \cite[Lemma 1.3]{I1} that 
each $\gamma \in \Gamma_{\Delta} - \{ 1 \}$ has its attractive (resp. repulsive) 
fixed points $\alpha$ (resp. $\alpha'$) in ${\mathbb P}^{1}_{B_{\Delta}}$ 
and its multiplier $\beta \in \sum_{e \in E} A_{\Delta} \cdot y_{e}$ which satisfy 
$$
\frac{\gamma(z) - \alpha}{z - \alpha} = \beta \frac{\gamma(z) - \alpha'}{z - \alpha'}. 
$$ 
Fix a set $\{ \gamma_{1},..., \gamma_{g} \}$ of generators of $\Gamma_{\Delta}$, 
and for each $\gamma_{i}$, 
denote by $\alpha_{i}$ (resp. $\alpha_{-i}$) its attractive (resp. repulsive) fixed points, 
and by $\beta_{i}$ its multiplier.  
Then under the assumption that there is no element of $\pm E \cap {\cal E}_{\infty}$ 
with terminal vertex $v_{b}$, 
for each $1 \leq i \leq g$, 
we define the associated {\it universal differential of the first kind} as 
$$
\omega_{i} = 
\sum_{\gamma \in \Gamma_{\Delta} / \left\langle \gamma_{i} \right\rangle} 
\left( \frac{1}{z - \gamma(\alpha_{i})} - \frac{1}{z - \gamma(\alpha_{-i})} \right) dz. 
$$
Assume that 
$$
\{ h \in \pm E \cap {\cal E}_{\infty} \ | \ v_{h} = v_{b} \} = \emptyset,  \ \ 
\{ t \in T \ | \ v_{t} = v_{b} \} \neq \emptyset. 
$$ 
Then for each $t \in T$ with $v_{t} = v_{b}$ and $k > 1$, 
we define the associated {\it universal differential of the second kind} as 
$$
\omega_{t, k} = \left\{ \begin{array}{ll} 
{\displaystyle \sum_{\gamma \in \Gamma_{\Delta}} \frac{d \gamma(z)}{(\gamma(z) - x_{t})^{k}}} 
& (x_{t} \neq \infty), 
\\
\\
{\displaystyle \sum_{\gamma \in \Gamma_{\Delta}} \gamma(z)^{k-2} d \gamma(z)} 
& (x_{t} = \infty). 
\end{array} \right. 
$$
Furthermore, 
put $T_{\infty} = \{ t \in T \cap {\cal E}_{\infty} \ | \ v_{t} = v_{b} \}$ 
whose cardinality is $0$ or $1$, 
and take a maximal subtree ${\cal T}_{\Delta}$ of $\Delta$, 
and for each $t \in T$, 
take the unique path $\rho_{t} = h(1) \cdots h(l)$ in ${\cal T}_{\Delta}$ 
from $v_{t}$ to $v_{b}$, 
and put $\phi_{t} = \phi_{h(l)} \cdots \phi_{h(1)}$. 
Then for each $t_{1}, t_{2} \in T$ with $t_{1} \neq t_{2}$, 
we define the associated {\it universal differential of the the third kind} as 
$$
\omega_{t_{1}, t_{2}} = \sum_{\gamma \in \Gamma_{\Delta}} 
\left( \frac{d \gamma(z)}{\gamma(z) - \phi_{t_{1}}(x_{t_{1}})} - 
\frac{d \gamma(z)}{\gamma(z) - \phi_{t_{2}}(x_{t_{2}})} \right), 
$$
where $\phi_{t_{i}}(x_{t_{i}}) = \infty$ if $t_{i} \in T_{\infty}$. 
\vspace{2ex}

\noindent
{\bf Theorem 3.1.} 
\begin{it} 

{\rm (1)} 
For each $1 \leq i \leq g$, 
$\omega_{i}$ is a regular differential on $C_{\Delta} \otimes_{A_{\Delta}} B_{\Delta}$ 
(cf. \cite[\S 3)]{MD}). 

{\rm (2)} 
For each $t \in T$ with $v_{t} = v_{b}$ and $k > 1$, 
$\omega_{t, k}$ is a meromorphic differential on $C_{\Delta} \otimes_{A_{\Delta}} B_{\Delta}$ 
which has only pole (of order $k$) at the point $p_{t}$ corresponding to $t$. 

{\rm (3)}
For each $t_{1}, t_{2} \in T$ such that $t_{1} \neq t_{2}$, 
$\omega_{t_{1}, t_{2}}$ is a meromorphic differential on 
$C_{\Delta} \otimes_{A_{\Delta}} B_{\Delta}$ 
which has only (simple) poles at the points $p_{t_{1}}$ (resp. $p_{t_{2}}$) 
corresponding to $t_{1}$ (resp $t_{2}$) with residue $1$ (resp. $-1$). 

{\rm (4)} 
Take $x_{h}$ $(h \in \pm E \cup T)$ and $y_{e}$ $(e \in E)$ 
be complex numbers as in 2.3. 
Then $\omega_{i}, \omega_{t, k}, \omega_{t_{1}, t_{2}}$ are abelian differentials 
on the Riemann surface $R_{\Gamma}$, 
where $\Gamma$ is the Schottky group obtained from $\Gamma_{\Delta}$. 

\end{it} 
\vspace{2ex}

\noindent
{\bf Proof.} 
By Proposition 2.1, 
$\omega_{i}$ are differentials on ${\cal P}_{T_{\Delta}}$, 
and for any $\delta \in \Gamma_{\Delta}$, 

\begin{eqnarray*}
\omega_{i}(\delta(z)) & = & 
\sum_{\gamma \in \Gamma_{\Delta} / \left\langle \gamma_{i} \right\rangle} 
\left( \frac{1}{\delta(z) - \gamma(\alpha_{i})} - \frac{1}{\delta(z) - \gamma(\alpha_{-i})} \right) 
d \delta(z) 
\\ 
& = & 
\sum_{\gamma \in \Gamma_{\Delta} / \left\langle \gamma_{i} \right\rangle} 
\left( \frac{\gamma(\alpha_{i}) - \gamma(\alpha_{-i})}
{(\delta(z) - \gamma(\alpha_{i})) (\delta(z) - \gamma(\alpha_{-i}))} \right) 
d \delta(z) 
\\ 
& = & 
\sum_{\gamma \in \Gamma_{\Delta} / \left\langle \gamma_{i} \right\rangle} 
\left( \frac{(\delta^{-1} \gamma)(\alpha_{i}) - (\delta^{-1} \gamma)(\alpha_{-i})}
{(z - (\delta^{-1} \gamma)(\alpha_{i})) (z - (\delta^{-1} \gamma)(\alpha_{-i}))} \right) dz 
\\ 
& = & 
\sum_{\gamma \in \Gamma_{\Delta} / \left\langle \gamma_{i} \right\rangle} 
\left( \frac{1}{z - (\delta^{-1} \gamma)(\alpha_{i})} - 
\frac{1}{z - (\delta^{-1} \gamma)(\alpha_{-i})} \right) dz 
\\ 
& = & 
\omega_{i}(z). 
\end{eqnarray*}
Therefore, by construction, 
$\omega_{i}$ are differentials on $C_{\Delta}$ 
which are regular outside the union of $y_{e} = 0$ $(e \in E)$, 
and hence the assertions (1) follows. 
One can prove (2)--(3) similarly, 
and we prove (4). 
As is stated in 2.1, 
$R_{\Gamma}$ is given by the quotient space $\Omega_{\Gamma}/\Gamma$. 
Under the assumption on complex numbers $x_{h}$ and $y_{e}$, 
it is shown in \cite{Scho} that $\sum_{\gamma \in \Gamma} |\gamma'(z)|$ 
is uniformly convergent on any compact subset in 
$\Omega_{\Gamma} - \cup_{\gamma \in \Gamma} \gamma(\infty)$, 
and hence the assertion holds for $\omega_{t, k}$. 
If $a \in \Omega_{\Gamma} - \cup_{\gamma \in \Gamma} \gamma(\infty)$, 
then $\lim_{n \rightarrow \infty} \gamma_{i}^{\pm n}(a) = \alpha_{\pm i}$, 
and hence 
\begin{eqnarray*}
d \left( \int_{a}^{\gamma_{i}(a)} \sum_{\gamma \in \Gamma} 
\frac{d \gamma(\zeta)}{\gamma(\zeta) - z} \right) 
& = & 
\sum_{\gamma \in \Gamma} 
\left( \frac{1}{z - (\gamma \gamma_{i})(a)} - \frac{1}{z - \gamma(a)} \right) dz 
\\ 
& = & 
\sum_{\gamma \in \Gamma / \langle \gamma_{i} \rangle} 
\lim_{n \rightarrow \infty} 
\left( \frac{1}{z - (\gamma \gamma_{i}^{n})(a)} - 
\frac{1}{z - (\gamma \gamma_{i}^{-n})(a)} \right) dz 
\\ 
& = & 
\omega_{i}(z). 
\end{eqnarray*}
Therefore, 
$\omega_{i}$ is absolutely and uniformly convergent 
on any compact subset in $\Omega_{\Gamma}$, 
and hence is an abelian differential on $\Omega_{\Gamma}/\Gamma$. 
\ $\square$ 
\vspace{2ex}

\noindent
{\it 3.2. Stability of universal differentials} 
\vspace{2ex}

For a vertex $v \in V$, 
denote by $C_{v}$ the corresponding irreducible component of 
$C_{\Delta} \otimes_{A_{\Delta}} A_{0}$. 
Then $P_{v} = {\mathbb P}^{1}_{A_{0}}$ is the normalization of $C_{v}$. 
\vspace{2ex}

\noindent
{\bf Theorem 3.2.} 
\begin{it} 

{\rm (1)} 
For each $1 \leq i \leq g$, 
let $\phi_{h_{i}(1)} \cdots \phi_{h_{i}(l_{i})}$ $(h_{i}(j) \in \pm E)$ be the unique reduced product 
such that $v_{-h_{i}(j)} = v_{h_{i}(j+1)}$ and $h_{i}(1) \neq - h_{i}(l_{i})$ which is conjugate to $\gamma_{i}$. 
Then for each $v \in V$, 
the pullback $\left( \omega_{i}|_{C_{v}} \right)^{*}$ of $\omega_{i}|_{C_{v}}$ to $P_{v}$ 
is given by 
$$
\left( \sum_{v_{h_{i}(j)} = v} \frac{1}{z - x_{h_{i}(j)}} - 
\sum_{v_{-h_{i}(k)} = v} \frac{1}{z - x_{-h_{i}(k)}} \right) dz. 
$$

{\rm (2)}
For each $v \in V$, 
$\left( \omega_{t, k}|_{C_{v}} \right)^{*}$ is given by 
$\displaystyle \frac{dz}{(z - x_{t})^{k}}$ if $v = v_{t}$, 
and is $0$ otherwise. 

{\rm (3)} 
Denote by $\rho_{t_{j}} = h_{j}(1) \cdots h_{j}(l_{j})$ 
the unique path from $v_{t_{j}}$ $(t_{j} \in T)$ to $v_{b}$ in ${\cal T}_{\Delta}$. 
Then for each $v \in V$, 
$\left( \omega_{t_{1}, t_{2}}|_{C_{v}} \right)^{*}$ is given by 
$$
\left( \sum_{v_{h} = v} \frac{1}{z - x_{h}} - 
\sum_{v_{-k} = v} \frac{1}{z - x_{-k}} \right) dz, 
$$
where $h, k$ runs through 
$\left\{ t_{1}, h_{1}(1),..., h_{1}(l_{1}), -t_{2}, -h_{2}(1),..., -h_{2}(l_{2}) \right\}$. 
\end{it} 
\vspace{2ex}

\noindent
{\bf Proof.} 
For the proof of (1), 
we may assume that 
$\gamma_{i} = \phi_{h_{i}(1)} \cdots \phi_{h_{i}(l_{i})}$. 
Let $\gamma$ be an element of $\Gamma_{\Delta}$. 
Then by Proposition 2.1 (2), putting $y_{e} = 0$ $(e \in E)$, 
$$
\frac{1}{z - \gamma(\alpha_{i})} - \frac{1}{z - \gamma(\alpha_{-i})} = 
\frac{\gamma(\alpha_{i}) - \gamma(\alpha_{-i})}
{(z - \gamma(\alpha_{i}))(z - \gamma(\alpha_{-i}))}  
$$
becomes 
$$
\frac{1}{z - x_{h_{i}(j)}} - \frac{1}{z - x_{-h_{i}(j-1)}}
$$ 
if $j \in \{ 1,..., l_{i} \}$ $\left( \mbox{$h_{i}(j-1) := h_{i}(l_{i})$ when $j = 1$} \right)$, 
$v(h_{i}(j)) = v$, $\phi_{h_{i}(j)} \phi_{h_{i}(j+1)} \cdots \phi_{h_{i}(l_{i})}$ belongs to 
$\gamma \langle \gamma_{i} \rangle$, 
and becomes $0$ otherwise. 
Therefore, the assertion follows from the definition of $\omega_{i}$. 

The assertion (2) follows from Proposition 2.1 (1) and (3), 
and the assertion (3) can be shown in the same way as above. 
\ $\square$ 
\vspace{2ex}

\noindent
{\bf Remark.} 
Let the notation be as above, 
and for a family of Schottky uniformized Riemann surfaces obtained from $C_{\Delta}$ 
as in 2.3, 
let $c_{e}$ be a cycle corresponding to $e \in E - E'$ 
which is oriented by the right-hand rule for $\gamma_{j}$. 
Then the restriction of $\omega_{j}$ to an irreducible component of $C_{\Delta'}$ 
is characterized analytically by that its integral along $c_{e}$ is $2 \pi \sqrt{-1}$ 
if $e = |h_{j,k}|$ for some integer $1 \leq k \leq l(j)$, 
and is $0$ otherwise. 
Therefore, Theorem 3.2 can be extended for this family 
degenerating as $y_{e} \rightarrow 0$ $(e \in E')$. 
Since Schottky uniformized Riemann surfaces make a nonempty open subset 
in the moduli space of Riemann surfaces, 
by the theorem of identity in the function theory of several variables, 
the extended version of Theorem 3.2 also holds 
for families of general Riemann surfaces. 
This modification gives a more explicit formula than \cite[Corollary 4.6]{HuN}. 
\vspace{2ex}

A {\it stable differential} on a stable curve is defined as (regular or meromorphic) 
section of the dualizing sheaf on this curve (\cite{DM}). 
By Theorems 3.1 and 3.2, we have: 
\vspace{2ex}

\noindent
{\bf Theorem 3.3.} 
\begin{it} 

{\rm (1)} 
For each $1 \leq i \leq g$, 
$\omega_{i}$ is a regular stable differential on $C_{\Delta}$. 

{\rm (2)} 
For each $t \in T$ with $v_{t} = v_{b}$ and $k > 1$, 
$\omega_{t, k}$ is a meromorphic stable differential on $C_{\Delta}$ 
which has only pole (of order $k$) at the point $p_{t}$ corresponding to $t$. 

{\rm (3)}
For each $t_{1}, t_{2} \in T$ with $t_{1} \neq t_{2}$, 
$\omega_{t_{1}, t_{2}}$ is a meromorphic stable differential on $C_{\Delta}$ 
which has only (simple) poles at the points $p_{t_{1}}$ (resp. $p_{t_{2}}$) 
corresponding to $t_{1}$ (resp $t_{2}$) with residue $1$ (resp. $-1$). 

\end{it} 
\vspace{4ex}

\noindent 
{\bf 4. Period map,  Gauss-Manin connection and variation of Hodge structure} 
\vspace{2ex}

\noindent
{\it 4.1. Period map} 
\vspace{2ex}

Let $\Delta = (V, E, T)$ be a stable graph, 
and assume that there are a vertex $v_{b} \in V$ satisfying 
$$
\{ h \in \pm E \cap {\cal E}_{\infty} \ | \ v_{h} = v_{b} \} = \emptyset,  \ \ 
\{ t \in T \ | \ v_{t} = v_{b} \} \neq \emptyset.  
$$
Then one can take an element of $T$ with terminal vertex $v_{b}$ 
which we denote by $t_{0}$. 
\vspace{2ex}

\noindent
{\bf Theorem 4.1.} 
\begin{it} 
Let $C_{\Delta}^{\circ} = C_{\Delta} - \{ p_{t} \ | \ t \in T \}$ be the open curve over $A_{\Delta}$. 
Then 
$$
\omega_{i} \ (1 \leq i \leq g), \ \ \omega_{t_{0}, k} \ (2 \leq k \leq g+1), \ \ 
\omega_{t, t_{0}} \ (t \in T - \{ t_{0} \})  
$$ 
make a basis of $H^{1}_{\rm dR} \left( C_{\Delta}^{\circ}/B_{\Delta}, K_{\Delta} \right)$, 
where $K_{\Delta}$ denotes the quotient field of $A \widehat{\otimes} {\mathbb Q}$. 
\end{it}
\vspace{2ex}

\noindent
{\bf Proof.} 
To prove this theorem, 
we can regard $C_{\Delta}^{\circ}$ as a family of open Riemann surfaces 
obtained from $C_{\Delta}$ by removing sections associated with $p_{t}$ $(t \in T)$.  
For small (counter-clockwise oriented) closed paths $a_{i}$ $(1 \leq i \leq g)$ 
in $C_{\Delta}^{\circ}$ around $\alpha_{i}$, 
$$
\int_{a_{i}} \omega_{j} = 2 \pi \sqrt{-1} \delta_{i, j}, 
\ \ 
\int_{a_{i}} \omega_{t_{0}, k} = \int_{a_{i}} \omega_{t, t_{0}} = 0, 
$$
where $\delta_{i, j}$ denotes the Kronecker delta. 
Furthermore, for small closed paths $a_{t}$ $(t \in T)$ in $C_{\Delta}^{\circ}$ 
around $p_{t}$, 
$$
\int_{a_{t}} \omega_{t', t_{0}} = 
2 \pi \sqrt{-1} \delta_{t, t'} = - 2 \pi \sqrt{-1} \delta_{t, t_{0}}, 
\ \ 
\int_{a_{t}} \omega_{j} = \int_{a_{t}} \omega_{t_{0}, k} = 0. 
$$
Let $b_{i}$ $(1 \leq i \leq g)$ be a closed path in $C_{\Delta}^{\circ}$ 
corresponding to a path from $p_{\infty}$ to $\gamma_{i}(p_{\infty})$. 
Then by Proposition 2.1 (1) and (3), 
${\displaystyle \int_{b_{i}} \omega_{t_{0}, k}}$ $(k > 1)$ 
belongs to $A_{\Delta} \widehat{\otimes} {\mathbb Q}$ whose constant term is 
$$
\int_{\infty}^{\alpha_{i}} \frac{dz}{(z - x_{t_{0}})^{k}} = 
\frac{1}{(-k+1) (\alpha_{i} - x_{t_{0}})^{k-1}}.
$$
We will prove that  
$$
\Omega := \det \left( \int_{b_{i}} \omega_{t_{0}, j+1} \right)_{1 \leq i, j \leq g} \neq 0. 
$$
When the set $V$ of vertices in $\Delta$ consists of one element, 
$\alpha_{i} - \alpha_{j} \in A_{\Delta}^{\times}$ for $i \neq j$ which implies that 
the Vandermonde determinant 
$$
\det \left( \frac{1}{(\alpha_{i} - x_{t_{0}})^{j}} \right)_{1 \leq i, j \leq g} = 
\prod_{i=1}^{g} \frac{1}{\alpha_{i} - x_{t_{0}}} 
\prod_{1 \leq i < j \leq g} \frac{\alpha_{i} - \alpha_{j}}{(\alpha_{i} - x_{t_{0}})(\alpha_{j} - x_{t_{0}})}
$$ 
belongs to $A_{\Delta}^{\times}$, and hence $\Omega \neq 0$. 
Denote by ${\cal M}_{g}$ the moduli stack of proper smooth curves of genus $g$ 
which is shown to be geometrically irreducible (cf. \cite{DM}). 
For any $\Delta$, 
the image of the morphism ${\rm Spec}(B_{\Delta}) \rightarrow {\cal M}_{g}$ 
associated with $C_{\Delta}/B_{\Delta}$ is Zariski dense, 
and hence the evaluation of $\Omega$ on ${\rm Spec}(B_{\Delta})$ is also not $0$. 
Therefore, 
$\left\{ \omega_{i}, \omega_{t_{0}, i+1} \right\}_{1 \leq i \leq g}$ gives a basis of 
$H^{1}_{\rm dR} (C_{\Delta}/B_{\Delta}, K_{\Delta})$ from which the assertion follows. 
\ $\square$ 
\vspace{2ex}

\noindent
{\it 4.2. Gauss-Manin connectiom and variation of Hodge structure} 
\vspace{2ex}

Let $\nabla$ denote the Gauss-Manin connection 
for $C_{\Delta}$ regarded as a family of proper smooth curves over $B_{\Delta}$. 
\vspace{2ex} 

\noindent
{\bf Theorem 4.2.} 
\begin{it} 

{\rm (1)} 
There exist elements $\eta_{j}$ $(1 \leq j \leq g)$ of 
$H^{1}_{\rm dR}  \left( C_{\Delta}/B_{\Delta}, K_{\Delta} \right)$ which are represented as 
$K_{\Delta}$-linear sums of $\omega_{t_{0}, k}$ $(2 \leq k \leq g+1)$, and satisfy 
$$
\nabla (\omega_{i}) = \sum_{j=1}^{g} \eta_{j} \otimes (d P_{ij}/P_{ij}), \ \nabla (\eta_{i}) = 0 
\ \ (1 \leq i \leq g),  
$$  
where $P_{ij} \in B_{\Delta}$ denote the universal periods \cite[Section 3]{I1} defined as 
$\exp \left( \int_{\gamma_{i}} \omega_{j} \right)$. 

{\rm (2)} 
The both sets $\{ \eta_{j} \ | \ 1 \leq j \leq g\}$ and 
$\{ \omega_{t_{0}, k} \ | \ 2 \leq k \leq g+1 \}$ give bases 
of the Hodge component $H^{0, 1}$ of 
$H^{1}_{\rm dR}  \left( C_{\Delta}/B_{\Delta}, K_{\Delta} \right)$. 
\end{it}
\vspace{2ex}

\noindent
{\bf Proof.} 
By Theorem 4.1 and its proof, 
there are $K_{\Delta}$-linear sums $\eta_{i}$ of $\omega_{t_{0}, k}$ 
such that $\int_{a_{i}} \eta_{j} = 0$ and $\int_{b_{i}} \eta_{j} = \delta_{i,j}$. 
Then the assertion (1) follows. 
Therefore, the Griffiths transversality implies that $\{ \eta_{j} \}$ generates 
a sub $K_{\Delta}$-module of $H^{0,1}$ with rank $g$ 
which is also generated by $\{ \omega_{t_{0}, k} \}$. 
Then the assertion (2) follows. 
\ $\square$ 
\vspace{2ex}

\noindent
{\bf Corollary 4.3.} 
\begin{it} 
Let $C/S$ be a family of Mumford curve over a $p$-adic field of characteristic $0$ 
obtained from $C_{\Delta}$ as in 2.3. 
Then $\{ \omega_{i} \}$ (resp. $\{ \eta_{j} \}$) become 
$\{ \alpha_{i} \}$ (resp. $\{ \beta_{j} \}$) given in \cite[Theorem 2]{G} 
which are bases of the Hodge components $H^{1,0}(C)$ (resp. $H^{0,1}(C)$) 
given in \cite[Theorems 2.8 and 2.9]{dS1}. 
\end{it}
\vspace{2ex}

\noindent
{\bf Proof.} 
The universal differential $\omega_{i}$ becomes a differential of the first kind on $C$ 
whose residue (cf. \cite{S}) for each edge in $E$ is constant, 
and $\eta_{j}$ becomes a differential of the second kind on $C$ whose Coleman integration 
(cf. \cite{C, CdS}) along $\gamma_{i}$ is $\delta_{i, j}$. 
Then the assertion follows from Theorems 2.8 and 2.9 of \cite{dS1}. 
\ $\square$ 
\vspace{2ex}

\noindent
{\it 4.3. Monodromy weight filtration} 
\vspace{2ex}

Let $C_{\Delta}$ be the generalized Tate curves associated with 
a stable graph $\Delta = (V, E, T)$, 
and put 
$H = H^{1}_{\rm dR} \left( C_{\Delta}/B_{\Delta}, K_{\Delta} \right)$. 
For a subset $E'$ of $E$, denote by 
$$
\{ 0 \} \subset F_{E'}^{0} H \subset F_{E'}^{1} H \subset H, \ \ 
N_{E'} : H / F_{E'}^{1} H \stackrel{\sim}{\rightarrow} F_{E'}^{0} H
$$
the monodromy weight filtration and operator with respect to the closed fiber of 
$C_{\Delta}$ obtained by $y_{e}$ $(e \in E')$. 
\vspace{2ex}

\noindent
{\bf Theorem 4.4.} 
\begin{it} 

{\rm (1)} 
For each $e \in E$, 
let ${\rm Res}_{e}$ be the linear map on $H$ defined by 
${\rm Res}_{e} \left( \omega_{t_{0}, k} \right) = 0$ and by 
$$
{\rm Res}_{e} \left( \omega_{i} \right) = 
\sharp \{ j \ | \ h_{i}(j) = e \} - \sharp \{ k \ | \ h_{i}(k) = -e \}, 
$$
where $\phi_{h_{i}(1)} \cdots \phi_{h_{i}(l)}$ $(h_{i}(j) \in \pm E)$ denotes the reduced product 
with $v_{-h_{i}(j)} = v_{h_{i}(j+1)}$ which is conjugate to $\gamma_{i}$ $(1 \leq i \leq g)$. 
Then we have 
$$
F_{E'}^{1} H = {\rm Ker} \left( \bigoplus_{e \in E'} {\rm Res}_{e} \right). 
$$ 

{\rm (2)} 
Let $\Delta'$ be the graph obtained from $\Delta$ 
by shrinking each edge in $E - E'$ to one point. 
Then $\{ {\rm Res}_{e} \}_{e \in E'}$ gives rise to an isomorphism from $H/F_{E'}^{1} H$ 
onto the group $H^{1} \left( \Delta', K_{\Delta} \right)$ 
of singular cohomology classes. 
Furthermore, the line integral in $C_{\Delta}$ corresponding to cycles in $\Delta'$ 
gives an isomorphism 
$$
F_{E'}^{0} H \cong H^{1} \left( \Delta', K_{\Delta} \right), 
$$
and $N_{E'}$ gives the identity map on $H^{1} \left( \Delta', K_{\Delta} \right)$ 
under the above isomorphisms.  
\end{it}
\vspace{2ex}

\noindent
{\bf Proof.} 
The assertion (1) follows from that ${\rm Res}_{e}$ gives the residue map 
for the family of Riemann surfaces associated with $C_{\Delta}$. 
Then by Theorem 4.2, 
$F_{E}^{0} H = F_{E}^{1} H$ has a basis $\{ \omega_{t_{0}, k} \}$, 
and by Theorem 4.1, the line integrals 
from $p_{\infty}$ to $\gamma_{i}(p_{\infty})$ $(1 \leq i \leq g)$ gives an isomorphism 
$F_{E}^{0} H \stackrel{\sim}{\rightarrow} H^{1} \left( \Delta, K_{\Delta} \right)$. 
Furthermore, Theorem 4.2 implies that $N_{E}$ is the identity map on 
$H^{1} \left( \Delta, K_{\Delta} \right)$. 
Under this isomorphism, 
$F_{E'}^{0} H$ corresponds to the image of 
$H^{1} \left( \Delta', K_{\Delta} \right) \hookrightarrow H^{1} \left( \Delta, K_{\Delta} \right)$ 
which implies (2). 
\ $\square$ 
\vspace{2ex}

\noindent
{\it 4.4. Coleman integration in families of Mumford curves} 
\vspace{2ex}

We consider the problem of constructing Coleman integrals in families 
raised by Besser \cite[1.6]{Be} in the case of Mumford curves.   
Let $K$ be a $p$-adic field, namely a local field of characteristic $0$ 
whose residue field is of characteristic $p > 0$. 
Denote by $| \cdot |$ a multiplicative valuation on $K$, 
and by ${\cal S}_{\Delta}$ the $K$-analytic subspace of $K^{\cal E} \times K^{E}$ 
consists of $(x_{h}, y_{e})_{h \in {\cal E}, e \in E}$ satisfying 
$$
|x_{h}| \leq 1, \ \ |x_{e} - x_{-e}| = 1, \ \ 
|x_{h} - x_{h'}| = 1 \ (\mbox{$h, h \in {\cal E}$ with $h \neq h'$ and $v_{h} = v_{h'}$}). 
$$
Then as is seen in 2.3, 
we have a family ${\cal C}_{\Delta}/{\cal S}_{\Delta}$ of Mumford curves over $K$ 
from the generalized Tate curve $C_{\Delta}$ by specializing parameters. 
Therefore, the universal differentials $\omega_{i}$ and $\omega_{t, k}$ in 3.1 give differentials 
on ${\cal C}_{\Delta}/{\cal S}_{\Delta}$ of the first and second kind respectively, 
and hence there exist their Coleman integrals given by \cite{CdS} 
under fixing a logarithm homomorphism $\log : K^{\times} \rightarrow K$. 
\vspace{2ex}

\noindent
{\bf Theorem 4.5.} 
\begin{it}
The Coleman integrals of $\omega_{i}$ and $\omega_{t, k}$ (determined unique up to constants)  
on ${\cal C}_{\Delta}/{\cal S}_{\Delta}$ are given by 
$$
\log \left( \prod_{\gamma \in \Gamma / \langle \gamma_{i} \rangle} 
\frac{z - \gamma(\alpha_{i})}{z - \gamma(\alpha_{-i})} \right) 
$$
and 
$$
\sum_{\gamma \in \Gamma} \frac{1}{1 - k} 
\left( \frac{1}{(\gamma(z) - x_{t})^{k-1}} - \frac{1}{(\gamma(z_{0}) - x_{t})^{k-1}} \right), 
$$ 
where $\Gamma$ denotes Schottky groups over $K$ associated with $\Gamma_{\Delta}$, 
and $z_{0}$ is a point on ${\mathbb P}^{1}_{K}$ outside the limit set of $\Gamma$. 
\end{it}
\vspace{2ex}

\noindent
{\bf Proof.} 
We use the description given by de Shalit \cite{dS2} of Coleman integrals 
on semistable curves. 
Let $C_{K}$ be a Mumford curve over $K$ as a member of ${\cal C}_{\Delta}$. 
Then it is shown in \cite[0.4 and 1.1]{dS2} that there exists a cover $\widetilde{C}_{K}$ 
of $C_{K}$ as a $K$-analytic space with action of $\Gamma$ such that 
$\widetilde{C}_{K}/\Gamma \cong C_{K}$, 
and the Coleman integrals on $C_{K}$ are functions on $\widetilde{C}_{K}$ 
obtained by anti-differentiating. 
Therefore, the assertion follows from the explicit formulas of $\omega_{i}$ and $\omega_{t, k}$. 
\ $\square$ 
\vspace{2ex}

\noindent
{\bf Corollary 4.6.} 
\begin{it} 
Let $C_{K}$ be a Mumford curve over a $p$-adic field $K$ obtained from $C_{\Delta}$
by the base change corresponding to a ring homomorphism $\varphi$ from 
$A_{\Delta}$ to the valuation ring of $K$, 
and put $E' = \{ e \in E \ | \ \varphi(y_{e}) = 0 \}$. 
Then the monodromy weight filtration and monodromy operator   
$$
\{ 0 \} \subset F_{E'}^{0} H \subset F_{E'}^{1} H \subset H, \ \ 
N_{E'} = H / F_{E'}^{1} H \rightarrow F_{E'}^{0} H. 
$$
given in Theorem 4.4 become the $p$-adic monodromy weight filtration and 
monodromy operator for $C_{K}$ (cf. \cite{CI}). 
\end{it}
\vspace{2ex}

\noindent
{\bf Proof.} The assertion follows from Theorem 4.5 and \cite[Theorem 0.1]{dS2}. 
\ $\square$ 
\vspace{4ex}

\noindent
{\bf 5. Unipotent period} 
\vspace{2ex}

\noindent
{\it 5.1. Unipotent fundamental group and period} 
\vspace{2ex}

Assume that $n$ is an integer $> 1$, 
and let $C^{\circ}$ be an algebraic curve over a subfield $K$ of ${\mathbb C}$ which is  obtained from a proper smooth curve $C$ of genus $g$ by removing $n$ points. 
Since $C^{\circ}$ is not complete, 
$H^{1}_{\rm dR}(C^{\circ})$ has a basis $B_{C^{\circ}}$ consisting of 
$2g + n-1$ meromorphic differentials on $C/K$ of the first or second kind 
which may have poles outside $C^{\circ}$. 
For $w_{1},..., w_{m} \in B_{C^{\circ}}$, 
we define a ${\cal D}$-module ${\cal D}(w_{1},..., w_{m})$ on $C^{\circ}$ 
whose underlying bundle is the trivial bundle $K^{m+1} \times C^{\circ}$ with connection form 
$d - \sum_{i = 1}^{m} e_{i, i+1} w_{i}$, 
where $e_{i,j}$ denotes the square matrix of degree $m+1$ whose $(k, l)$-entry is 
$\delta_{i,k} \cdot \delta_{j,l}$. 
We consider the tannakian subcategory of ${\cal D}$-modules on $C^{\circ}$ 
generated by ${\cal D}(w_{1},..., w_{m})$ $\left( w_{i} \in B_{C^{\circ}} \right)$. 
Since these underlying bundles are trivial, 
for each $K$-rational point $x$ on $C^{\circ}$, 
one can define the fiber functor on this category by taking the (trivial) fibers over $x$. 
Denote by $\pi_{1}^{\rm dR}(C^{\circ}; x)$ the tannakian fundamental group of 
this category which is a pro-finite algebraic group over $K$, 
and by ${\cal A}^{\rm dR}(C^{\circ}; x)$ the enveloping algebra of the Lie algebra 
${\rm Lie} \left( \pi_{1}^{\rm dR}(C^{\circ}; x) \right)$. 

Let $(C^{\circ})^{\rm an}$ be a Riemann surface associated with 
$C^{\circ} \otimes_{K} {\mathbb C}$, 
and for $K$-rational points $x, y$ on $C^{\circ}$, 
denote by $\pi_{1}((C^{\circ})^{\rm an}; x, y)$ the set of homotopy classes of paths 
from $x$ to $y$ in $(C^{\circ})^{\rm an}$. 
When $x = y$, $\pi_{1}((C^{\circ})^{\rm an}; x, y)$ becomes the fundamental group 
$\pi_{1}((C^{\circ})^{\rm an}; x)$ of $(C^{\circ})^{\rm an}$ with base point $x$. 
We consider the tannakian category of unipotent local systems on $(C^{\circ})^{\rm an}$ 
with fiber functor obtained from taking the fiber over $x$.  
Then it is shown in \cite{D2} that its tannakian fundamental group 
$\pi_{1}^{\rm Be}((C^{\circ})^{\rm an}; x)$ is a pro-algebraic group over ${\mathbb Q}$, 
and the associated enveloping algebra ${\cal A}^{\rm Be}((C^{\circ})^{\rm an}; x)$ 
of ${\rm Lie} \left( \pi_{1}^{\rm Be}((C^{\circ})^{\rm an}; x) \right)$ is isomorphic to 
$$
\lim_{m \rightarrow \infty} 
\left. {\mathbb Q} \left[ \pi_{1} ((C^{\circ})^{\rm an}; x) \right] \right/ I^{m},  
$$
where $I$ denotes the augmentation ideal of the group ring 
${\mathbb Q} \left[ \pi_{1}((C^{\circ})^{\rm an}; x) \right]$. 
Since $\pi_{1} ((C^{\circ})^{\rm an}; x)$ is a free group of rank $2g + n-1$, 
${\cal A}^{\rm Be} ((C^{\circ})^{\rm an}; x)$ becomes the ring of 
noncommutative formal power series over ${\mathbb Q}$ in $2g + n-1$ variables. 
For smooth $1$-forms $w_{1},..., w_{m}$ on $(C^{\circ})^{\rm an}$ 
and a smooth path $\gamma : (0, 1) \rightarrow (C^{\circ})^{\rm an}$, 
the iterated integral $\int_{\gamma} w_{1} \cdots w_{m}$ is defined as 
$$
\int_{0 < t_{1} < \cdots < t_{m} < 1} \gamma^{*}(w_{1})(t_{1}) \cdots \gamma^{*}(w_{m})(t_{m}). 
$$ 

\noindent
{\bf Theorem 5.1.} 
\begin{it}
There exists a canonical isomorphism 
\end{it}
$$
\pi_{1}^{\rm dR} (C^{\circ}; x, y) \otimes_{K} {\mathbb C} \cong 
\pi_{1}^{\rm Be} ((C^{\circ})^{\rm an}; x, y) \otimes_{\mathbb Q} {\mathbb C} 
$$
expressed by iterated integrals. 
\vspace{2ex} 

\noindent
{\bf Proof.} 
We may show the assertion in the case when $x = y$. 
For each ${\cal D}$-module ${\cal D}(w_{1},..., w_{m})$, 
there exists the associate local system on $(C^{\circ})^{\rm an}$ described by 
the iterated integrals $\int w'_{1} \cdots w'_{r}$ for $r \leq m$ and 
$w'_{i} \in \{ w_{1},..., w_{m} \}$. 
By this association, 
one has a group homomorphism 
$$
\pi_{1}^{\rm Be} ((C^{\circ})^{\rm an}; x) \otimes_{\mathbb Q} {\mathbb C} 
\rightarrow 
\pi_{1}^{\rm dR} (C^{\circ}; x) \otimes_{K} {\mathbb C}  
$$
which gives a ring homomorphism 
$$
\varphi : {\cal A}^{\rm Be} ((C^{\circ})^{\rm an}; x) \otimes_{\mathbb Q} {\mathbb C} 
\rightarrow {\cal A}^{\rm dR} (C^{\circ}; x) \otimes_{K} {\mathbb C}.
$$
Let $B_{m}((C^{\circ})^{\rm an})$ be the ${\mathbb C}$-vector space of iterated integrals spanned by 
$$
\int w_{1} \cdots w_{r} \ 
\left( w_{i} \in B_{C^{\circ}}, \ r \leq m \right),  
$$
and $H^{0} \left( B_{m}((C^{\circ})^{\rm an}); x \right)$ be the space consisting of 
elements of $B_{m}((C^{\circ})^{\rm an})$ whose restriction to 
$\left\{ \mbox{loops in $(C^{\circ})^{\rm an}$ based at $x$} \right\}$ is homotopy functional. 
Then it is shown in \cite[(5.3)]{H1} that 
$H^{0} \left( B_{m}((C^{\circ})^{\rm an}); x \right)$ is the dual space of 
$\left. {\mathbb C} \left[ \pi_{1} ((C^{\circ})^{\rm an}; x) \right] \right/ I^{m+1}$, 
and that by taking leading terms of iterated integrals, 
one has an exact sequence 
$$
0 \rightarrow H^{0} \left( B_{m-1}((C^{\circ})^{\rm an}); x \right) \rightarrow 
H^{0} \left( B_{m}((C^{\circ})^{\rm an}); x \right) \rightarrow 
H^{1}_{\rm dR}((C^{\circ})^{\rm an})^{\otimes m}. 
$$ 
Denote by $(I^{m}/I^{m+1})^{\vee}$ the dual space of $I^{m}/I^{m+1}$. 
Then there exists the associated injective linear map 
$$
(I^{m}/I^{m+1})^{\vee} \otimes_{\mathbb Q} {\mathbb C} \rightarrow 
H^{1}_{\rm dR}((C^{\circ})^{\rm an})^{\otimes m} \cong 
H^{1}_{\rm dR}(C^{\circ})^{\otimes m} \otimes_{K} {\mathbb C} 
$$ 
which is the inverse of the map 
$$
H^{1}_{\rm dR}(C^{\circ})^{\otimes m} \otimes_{K} {\mathbb C} \rightarrow (I^{m}/I^{m+1})^{\vee} \otimes_{\mathbb Q} {\mathbb C} 
$$ 
induced from the above homomorphism $\varphi$. 
Therefore,  $\varphi$ is an isomorphism. 
The remaining assertions follow from that the tannakian fundamental group $\pi_{1}^{\sharp}$ 
$(\sharp = {\rm dR}, {\rm Be})$ and its Lie algebra are given as the subsets of 
${\cal A}^{\sharp}$ which consist of grouplike elements and primitive elements in 
${\cal A}^{\sharp}$ respectively. 
\ $\square$ 
\vspace{2ex}

Following \cite{H1}, 
we define the Hodge (resp. weight) filtrations $F^{\bullet}$ (resp. $W_{\bullet}$) on 
$\pi^{\rm dR}_{1} (C^{\circ}; x, y)$ as follows. 
First, $F^{p}$ is spanned by ${\cal D}$-modules ${\cal D}(w_{1},..., w_{r})$  
for at least $p$ elements of $B_{C^{\circ}}$ which are holomorphic on $C$ or 
having only simple poles at $D = C - C^{\circ}$. 
Second, $W_{l}$ is spanned by ${\cal D}(w_{1},..., w_{r})$ for at most $l$ elements of $B_{C^{\circ}}$  which have only simple poles at $D$. 

Take a point $x$ on $C$, 
and a set $W$ of $K$-basis of $H_{\rm dR}^{1}(C^{\circ})$ consisting of 
meromorphic differentials on $C$ which are regular or having a simple pole at $x$. 
We associate the symbols $A_{w}$ with $w \in W$, 
and consider the generalized Knizhnik-Zamolodchikov (KZ for short) equation 
\begin{eqnarray*}
d F(z) = \left( \sum_{w \in W} A_{w} \cdot w \right) F(z). 
\end{eqnarray*}
Denote by $W_{x}$ the set of elements of $W$ with simple pole at $x$. 
Then the solution $F(z)$ of the generalized KZ equation normalized as 
$$
\lim_{(z-x) \downarrow 0} F(z) \cdot 
(z-x)^{- \sum_{w \in W_{x}} {\rm Res}_{x}(w) A_{w}} = 1 
$$ 
is described by a non-commutative formal power series in $A_{w}$ $(w \in W)$ 
whose coefficients are iterated integrals 
$$
\int_{x}^{z} w_{1} \cdots w_{m} \ \ (w_{i} \in W, \ w_{1} \not\in W_{x}) 
$$
which are called {\it unipotent periods} for $x$ with $W_{x} = \emptyset$ 
or the associated tangential point at $x$ with $W_{x} \neq \emptyset$. 
\vspace{2ex}

\noindent
{\it 5.2. Universal unipotent periods} 
\vspace{2ex}

For positive integers $k_{1},..., k_{l}$, we review the multiple polylogarithm function 
$$
{\rm Li}_{k_{1},..., k_{l}}(z) = 
\sum_{0 < n_{1} < \cdots < n_{l}} \frac{z^{n_{l}}}{n_{1}^{k_{1}} \cdots n_{l}^{k_{l}}} \ \ (|z| < 1) 
$$
which can be analytically continued to a multi-valued function 
on ${\mathbb P}^{1}_{\mathbb C} - \{ 0, 1, \infty \}$ 
which is represented as the iterated integral 
$$
{\rm Li}_{k_{1},..., k_{l}}(z) = 
\int_{0}^{z} w_{1} \underbrace{w_{0} \cdots w_{0}}_{k_{1}-1} 
w_{1} \underbrace{w_{0} \cdots w_{0}}_{k_{2}-1} w_{1} \cdots w_{1} 
\underbrace{w_{0} \cdots w_{0}}_{k_{l}-1}, 
$$
where $w_{0} = dz/z$, $w_{1} = dz/(1-z)$. 
Therefore, as is shown by Deligne \cite{D2} and Drinfe'ld \cite{Dr}, 
multiple polylogarithm functions give rise to local unipotent periods for ${\mathbb P}^{1} - \{ 0, 1, \infty \}$ 
which are considered as period functions of mixed Tate motives over ${\mathbb Z}$. 
\vspace{2ex}

\noindent
{\bf Proposition 5.2.} 
\begin{it} 
Let $w_{1},..., w_{m}$ be elements of 
$$ 
\left\{ \left. z^{k} dz, \ (z - 1)^{-l} dz \ \right| 
\ k, l \in {\mathbb Z}, \ k \geq -1, \ l \geq 1 \right\}. 
$$
Then the iterated integral 
$$
\int_{0}^{z} w_{1} \cdots w_{m}  \ \left( w_{1} \neq dz/z \right) 
$$
becomes a linear sum over the ring $R = {\mathbb Q} [z, 1/(z-1)]$ 
of multiple polylogarithm functions ${\rm Li}_{k_{1},..., k_{l}}(z)$. 
\end{it}  
\vspace{2ex} 

\noindent
{\bf Proof.} 
Since the derivative of ${\rm Li}_{k_{1},..., k_{l}}(z)$ is given by 
$$
\left\{ \begin{array}{ll} 
{\displaystyle \frac{{\rm Li}_{k_{1},..., k_{l}-1}(z)}{z}} & (k_{l} > 1), 
\\
{\displaystyle \frac{{\rm Li}_{k_{1},..., k_{l-1}}(z)}{1-z}} & (k_{l} = 1), 
\end{array} \right. 
$$
the assertion is proved by induction on $m$ using integration by parts. 
\ $\square$  
\vspace{2ex}

Let the notation be as in 4.1, 
and assume that $\Delta$ is trivalent, 
namely $C_{\Delta} \otimes_{A} A_{0}$ is maximally degenerate. 
Then for each $v \in V$, 
one can take 
$$
\{ x_{h} \ | \ h \in \pm E \cup T, \ v_{h} = v \} = \{ 0, 1, \infty \}, 
$$
hence $A_{0} = {\mathbb Z}$, 
and $B_{\Delta} \widehat{\otimes} {\mathbb Q}$ is the ring ${\mathbb Q}((y_{e}))$ 
of Laurent power series in $y_{e}$ $(e \in E)$. 
For each $t \in T$, 
denote by $p_{t}$ the associated point on $C_{\Delta}$ from which 
there exist two tangential points over ${\mathbb Z}$. 
Let 
$$
\omega_{i} \ (1 \leq i \leq g), \ \ 
\omega_{t_{0}, k} \ (1 \leq k \leq g), \ \ 
\omega_{t, t_{0}} \ (t \in T - \{ t_{0} \})  
$$ 
be the universal differentials given in 3.1 which make a basis of  
$H_{\rm dR}^{1} \left( C_{\Delta}/B_{\Delta} \widehat{\otimes} {\mathbb Q} \right)$
by Theorem 4.1. 

Assume that there is a tail in $T - \{ t_{0} \}$ which we denote by $t_{1}$, 
and take a tangential point $\vec{p}_{t_{1}}$ over ${\mathbb Z}$. 
Then the universal unipotent periods for $C_{\Delta}^{\circ}$ from $\vec{p}_{t_{1}}$ to $z$ 
are given by iterated integrals 
$$
\int_{p_{t_{1}}}^{z} w_{1} \cdots w_{m}, 
$$
where $w_{1},..., w_{m} \in \{ \omega_{i}, \omega_{t_{0}, k}, \omega_{t, t_{0}} \}$ 
such that $w_{1} \neq \omega_{t_{1}, t_{0}}$. 
\vspace{2ex}    

\noindent
{\bf Theorem 5.3.} 
\begin{it}
Take a system of above coordinates on $P_{v} = {\mathbb P}^{1}_{\mathbb C}$ $(v \in V)$ 
satisfying that $p_{t_{1}}$ is defined as $z = 0$ on $P_{v_{t_{1}}}$. 
Then the universal unipotent periods around $\vec{p}_{t_{1}}$ are represented as 
power series in $y_{e}$ $(e \in E)$ whose coefficients are 
${\mathbb Q}[z, 1/(z-1)]$-linear sums of multiple polylogarithm functions 
${\rm Li}_{k_{1},..., k_{l}}(z)$. 
\end{it}
\vspace{2ex}

\noindent
{\bf Proof.} 
By Proposition 2.1, 
$\omega_{i}$, $\omega_{t_{0}, k}$, $\omega_{t, t_{0}} (t \neq t_{1})$ and $\omega_{t_{1}, t_{0}} - dz/z$ 
are expanded as power series in $y_{e}$ whose coefficients are regular at $z = 0$, 
and hence the assertion follows from Propositions 2.1 (1) and 5.2.  
\ $\square$ 
\vspace{2ex}

\noindent
{\it 5.3. Universal $p$-adic unipotent periods} 
\vspace{2ex}

We review a result of Furusho \cite{F1} on the $p$-adic theory of multiple polylogarithm functions. 
Let $p$ be a prime number $p$, and denote by ${\mathbb C}_{p}$ the completion of 
the algebraic closure of the field ${\mathbb Q}_{p}$ of $p$-adic numbers. 
The field ${\mathbb C}_{p}$ has the complete valuation $| \ |_{p}$ normalized as $|p|_{p} = 1/p$. 
Fix an element $a$ of ${\mathbb C}_{p}$ which gives rise to a branch of $p$-adic logarithm 
$\log^{a} : {\mathbb C}_{p}^{\times} \rightarrow {\mathbb C}_{p}$ characterized by 
$\log^{a}(p) = a$. 
Then in \cite[Section 2]{F1}, the $p$-adic multiple polylogarithm function 
${\rm Li}^{a}_{k_{1},..., k_{l}}(z)$ is given as a locally analytic function on 
${\mathbb P}^{1}_{{\mathbb C}_{p}} - \{ 1, \infty \}$ satisfying that 
$$
{\rm Li}^{a}_{k_{1},..., k_{l}}(z) = 
\sum_{0 < n_{1} < \cdots < n_{l}} \frac{z^{n_{l}}}{n_{1}^{k_{1}} \cdots n_{l}^{k_{l}}} \ \ (|z|_{p} < 1). 
$$
It is shown by Furusho \cite{F1, F2} that local $p$-adic unipotent periods 
for ${\mathbb P}^{1} - \{ 0, 1, \infty \}$ are given by $p$-adic multiple polylogarithm functions. 

In the same way to \cite[Section 3]{F2}, 
we define the universal $p$-adic associator of the generalized KZ equation for 
$$
W = \left\{ \omega_{i} \ (1 \leq i \leq g), \ \omega_{t_{0}, k} \ (2 \leq k \leq g+1), \ 
\omega_{t, t_{0}} \ (t \in T - \{ t_{0} \}) \right\} 
$$
as the monodromy of normalized $p$-adic analytic solutions, 
and call the {\it universal $p$-adic unipotent periods} 
the coefficients of the universal $p$-adic associator. 
Then in the same way to showing Theorems 5.3 and 5.4, 
one can prove the following theorems which are examples of iterated Coleman integrals 
in families of Mumford curves (cf. \cite{Be}). 
\vspace{2ex}

\noindent
{\bf Theorem 5.4.} 
\begin{it}
Assume that there is a tail in $T - \{ t_{0} \}$ which we denote by $t_{1}$, 
and take a tangential point $\vec{p}_{t_{1}}$ over ${\mathbb Z}$. 
Take a system of above coordinates on $P_{v} = {\mathbb P}^{1}_{\mathbb C}$ $(v \in V)$ 
satisfying that $p_{t_{1}}$ is defined as $z = 0$. 
Then the $p$-adic unipotent periods around $\vec{p}_{t_{1}}$ are represented as power series in $y_{e}$ $(e \in E)$ whose coefficients are $p$-adic multiple polylogarithm functions 
${\rm Li}^{a}_{k_{1},..., k_{l}}(z)$. 
\end{it}

\vspace{4ex}

\noindent 
{\bf Acknowledgments} 
\vspace{2ex}

This work is partially supported by the JSPS Grant-in-Aid for 
Scientific Research No. 17K05179.

\renewcommand{\refname}{\centerline{\normalsize{\bf References}}}
\bibliographystyle{amsplain}

\end{document}